\documentclass[12pt]{amsart}
\usepackage[T1]{fontenc}
\usepackage{dsfont}
\usepackage{natbib}

\newtheorem{prop}{Proposition}
\newtheorem{lem}{Lemma}
\newtheorem{thm}{Theorem}
\newtheorem{cor}{Corollary}

\def\E{\mathbb{E}}
\def\var{{\mathrm{Var}}}
\def\S{\mathbb{S}}
\def\pen{{\mathrm{pen}}}
\def\PEN{{\mathrm{Pen}}}
\def\1{\mathds{1}}
\def\u{v}
\def\v{V}
\begin{document}

\title[]{Adaptive estimation of the transition density of a particular hidden Markov chain}

\maketitle

\begin{center}
\author{\textsc{Claire Lacour}\\
\vspace{0.5cm}
{MAP5, Universit\'e Paris 5, 45 rue des 
Saints-P\`eres\\ 75270 Paris Cedex 06, France.\\
email: lacour@math-info.univ-paris5.fr}
}\end{center}

\begin{abstract} We study the following model of hidden Markov chain: $Y_i=X_i+\varepsilon_i$, $ i=1,\dots,n+1$
with $(X_i)$ a real-valued positive recurrent and stationary Markov chain and $(\varepsilon_i)_{1\leq i\leq n+1}$ a noise independent of 
the sequence $(X_i)$ having a known distribution. We present an adaptive estimator of the transition density based on the quotient of 
a deconvolution estimator of the density of $X_i$ and an estimator of the density of $(X_i,X_{i+1})$. These estimators
are obtained by contrast minimization and model selection. We evaluate the $L^2$ risk and 
its rate of convergence for ordinary smooth and supersmooth noise with regard to ordinary smooth and supersmooth chains.
Some examples are also detailed.

\vspace{0.5cm}\noindent
\end{abstract}


\section{Introduction}

Let us consider the following model:
\begin{equation}
  Y_i=X_i+\varepsilon_i\qquad i=1,\dots,n+1
\end{equation}
where $(X_i)_{i\geq 1}$ is an irreducible and positive recurrent Markov chain and 
$(\varepsilon_i)_{i\geq 1}$ is a noise independent of $(X_i)_{i\geq 1}$. 
We assume that $\varepsilon_1,\dots, \varepsilon_n$ are independent and identically 
distributed random variables with known distribution. 
This model belongs to the class of hidden Markov models. Contrary to the literature on the subject,
we are interested in a nonparametric approach of the estimation of the hidden chain transition.  
The problem of estimating the density of $X_i$ from the observations 
$Y_1,\dots, Y_n$ when the $X_i$ are i.i.d. (known as the convolution model)
has been extensively studied, see e.g.
\cite{carrollhall}, \cite{fan91}, \cite{stefanski90}, \cite{penskyvidakovic}, \cite{comterota2006}.

But very few authors study the case where $(X_i)$ is a Markov chain. We can cite 
\cite{doreazhao02} who estimate the density of $Y_i$ in such a model, \cite{masry93}
who is interested in the estimation of the multivariate density in a mixing framework
and \cite{clemhmm}
who estimates the stationary density and the transition density of the hidden chain.
More precisely he introduces an estimator of the transition density 
based on thresholding of a wavelet-vaguelette decomposition and 
he studies its performance in the case of an ordinary smooth noise (i.e. whose Fourier 
transform has polynomial decay).
Here we are interested also in the estimation of the transition density of $(X_i)$
but we consider a larger class of noise distributions.
In \cite{clemhmm} there is no study of supersmooth noise (i.e. with exponentially decreasing Fourier transform), 
as the Gaussian distribution. However the study of such noise allows to find interesting rates 
of convergence, in particular  when the chain density is also supersmooth.
In the present paper, the four cases (ordinary smooth or supersmooth noise with 
ordinary smooth or supersmooth chain) are considered.

The aim of this paper is to estimate the transition density $\Pi$ of the Markov chain $(X_i)$
from the observations $Y_1,\dots, Y_n$. 
To do this we  assume that the regime is stationary and we note that $\Pi=F/f$ where $F$ is the density of $(X_i,X_{i+1})$
and $f$ the stationary density.
The estimation of $f$ comes down to a problem of deconvolution, 
as does the estimation of $F$. We use contrast minimization and a model selection method inspired by \cite{BBM} to find 
adaptive estimators of $f$ and $F$. Our estimator of $\Pi$ is then the quotient of the two 
previous estimators. Note that it is worth finding an adaptive estimator, i.e. an  estimator whose risk 
automatically achieves the minimax rates, because the regularity of the densities $f$ and $F$ is
generally very hard to compute, even if the chain can be fully described (case of a diffusion or an autoregressive process).

We study the performance of our estimator by computing the rate of convergence of the $L^2$ risk.
We improve the result of \cite{clemhmm} (case of an ordinary smooth noise) since we obtain the minimax rate
without logarithmic loss. Moreover we observe noteworthy rates of convergence in the case where
both noise and the chain are supersmooth.

The paper is organized as follows.
Section 2 is devoted to notations and assumptions while the estimation procedure is developed 
in Section 3. After describing the projection spaces to which the estimators belong, 
we define separately the estimator of the stationary density $f$, the one of the joined density
$F$ and last the estimator $\tilde{\Pi}$ of the transition density. Section 4 states the results
obtained for our estimators. 
To illustrate the theorems, some examples are provided in Section~5 as the AR(1) model, the Cox-Ingersoll-Ross process
or the stochastic volatility model. The proofs are to be found in Section 6.

\section{Notations and Assumptions}

For the sake of clarity, we use lowercase letters for the dimension 1 and capital letters for 
the dimension 2.
For a function $t: \mathbb{R}\mapsto \mathbb{R}$, we denote by $\|t\|$ the $L^2$ norm  that is 
$\|t\|^2=\int_\mathbb{R} t^2(x)dx.$
The Fourier transform $t^*$ of $t$ is defined by
$$t^*(u)=\int e^{-ixu} t(x)dx$$
Notice that the function $t$ is the inverse Fourier transform of $t^*$ and 
can be written $t(x)=1/(2\pi)\int e^{ixu} t^*(u)du$.
Finally, the convolution product is defined by 
${(t*s)(x)}=\int t(x-y)s(y)dy.$

In the same way, for a function $T:\mathbb{R}^2\mapsto \mathbb{R}$,
$\|T\|^2=\iint_{\mathbb{R}^2} T^2(x,y)dxdy$ and 
$$T^*(u,v)=\int\hspace{-3mm}\int e^{-ixu-iyv} T(x,y)dxdy,\quad
(T*S)(x,y)=\int\hspace{-3mm}\int T(x-z,y-w )S(z,w)dzdw.$$
We denote by $t\otimes s$ the function: $(x,y)\mapsto (t\otimes s)(x,y)=t(x)s(y)$.

The density of $\varepsilon_i$ is named $q$ and is considered as known. 
We denote by $p$ the density of $Y_i$. We have $p=f*q$ and then $p^*=f^*q^*$.
Similarly if $P$ is the density of $(Y_i, Y_{i+1 })$, then 
$P=F*(q\otimes q)$ and $P^*(u,v)=F^*(u,v)q^*(u)q^*(v)$.

Now the assumptions on the model are the following:
 \begin{description}
\item[A1] The function $q^*$ never vanishes. 
\item[A2] There exist $s\geq0, b> 0$, $\gamma>0$ if $s=0$ and $k_0, k_1>0$ such that
$$k_0(x^2+1)^{-\gamma/2}\exp(-b|x|^s)\leq |q^*(x)|\leq k_1(x^2+1)^{-\gamma/2}\exp(-b|x|^s)$$
\item[A3] The chain is stationary with (unknown) density $f$.
\item[A4] The chain is geometrically $\beta$-mixing ($\beta_q\leq M e^{-\theta q}$),
or arithmetically $\beta$-mixing ($\beta_q\leq M q^{-\theta}$) with $\theta>6$.
\end{description}
This last condition is verified as soon as the chain is uniformly ergodic.
In the sequel we consider the following smoothness spaces: 
$$\mathcal{A}_{\delta,r,a}(l)=\{f \text{ density on } \mathbb{R} \text{ and } 
\int |f^*(x)|^2(x^2+1)^{\delta}\exp(2a|x|^r)dx\leq l\}$$
with $r\geq0, a> 0, \delta>1/2$ if $r=0$, $l>0$ and 
\begin{eqnarray*}
& &  \mathbb{A}_{\Delta,R,A}(L)=\{F\text{ density on } \mathbb{R}^2 \text{ and } \\
& &\iint |F^*(x,y)|^2(x^2+1)^{\Delta}(y^2+1)^{\Delta}\exp(2A(|x|^R+|y|^R))dxdy\leq L\}
\end{eqnarray*}
with $R\geq0, A> 0, \Delta> 1/2$ if $R=0$, $L>0$.

When $r>0$ (respectively $R>0$) the function $f$ (resp. $F$) is known as supersmooth, and as ordinary smooth 
otherwise. In the same way, the noise distribution is called ordinary smooth if $s=0$ and supersmooth otherwise.
The spaces of ordinary smooth functions correspond to classic Sobolev classes, while
supersmooth functions are infinitely differentiable. It
includes for example normal ($r=2$) and Cauchy ($r=1$) densities.

\section{Estimation procedure}

Since $\Pi=F/f$ we proceed in 3 steps to estimate the transition density $\Pi$. 
First we find an estimator $\tilde{f}$ of $f$ (see
Subsection \ref{estimationf}). Then we estimate $F$ by $\tilde{F}$ (see
Subsection \ref{estimationF}). And last we estimate $\Pi$ with the quotient $\tilde{F}/\tilde{f}$
(Subsection \ref{estimationpi}).

All estimators defined here are projection estimators. We therefore start with describing the projection spaces.

\subsection{Projection spaces}
Let $$\varphi(x)=\sin(\pi x)/ (\pi x)$$ and, for $m$ in $\mathbb{N}^*$, $j$ in $\mathbb{Z}$,
 $\varphi_{m ,j}(x)=\sqrt{{m}}\varphi({m }x-j)$. 
Notice that $\{\varphi_{m,j}\}_{j\in \mathbb{Z}}$ is an orthonormal basis of the space of 
integrable functions having a Fourier transform with compact support included into
$[-\pi {m},\pi {m}].$
In the sequel, we use the following notations: 
$$S_m=\text{Span}\{\varphi_{m,j}\}_{j\in \mathbb{Z}}; \quad 
\mathbb{S}_m=\text{Span}\{\varphi_{m,j}\otimes\varphi_{m,k}\}_{j,k\in \mathbb{Z}}$$
These spaces have particular properties, which are a consequence of the first point of Lemma \ref{ordrevar}
(see Section \ref{techlem}):
\begin{eqnarray}\label{connex}
\forall t\in S_m \qquad \|t\|_\infty\leq \sqrt{{m}}\|t\|; \quad
\forall T\in \mathbb{S}_m \qquad \|T\|_\infty\leq {m}\|T\|
\end{eqnarray}
where $ \|t\|_\infty=\sup_{x\in\mathbb{R}}|t(x)|$
and $ \|T\|_\infty=\sup_{(x,y)\in\mathbb{R}^2}|T(x,y)|$.

\subsection{Estimation of $f$}\label{estimationf}

Here we estimate $f$, which is the density of the $X_i$'s. It is the classic deconvolution problem.
We choose to estimate $f$ by minimizing a contrast. The classical contrast in 
density estimation is 
$1/n \sum_{i=1}^n [\|t\|^2-2t(X_i)]$. It is not possible to use this contrast here
since we do not observe $X_1,\dots, X_n$. Only the noisy data 
$Y_1,\dots, Y_n$ are available. That is why we use the following lemma.

\begin{lem}\label{ut} For all function $t$, let $\u_t$ be the inverse Fourier transform of 
$t^*/q^*(-.)$, i.e. $$\u_t(x)=\frac{1}{2\pi}\int e^{ixu} \frac{t^*(u)}{q^*(-u)}du.$$ Then, for all $1\leq k\leq n$, 
  \begin{enumerate}
  \item $\E[\u_t(Y_k)|X_1,...,X_n]=t(X_k)$
  \item $\E[\u_t(Y_k)]=\E[t(X_k)]$
  \end{enumerate}
\end{lem}

The second assertion in Lemma \ref{ut} is an obvious consequence of the first one
and leads us to consider the following contrast:
\begin{equation*}
  \gamma_n(t)=\frac{1}{n} \sum_{i=1}^n [\|t\|^2-2\u_t(Y_i)]\quad
\text{ with }\quad \u_t^*(u)=\frac{t^*(u)}{q^*(-u)}
\end{equation*}
We can observe that $\E\gamma_n(t)={1}/{n} \sum_{i=1}^n[\|t\|^2-2\E[\u_t(Y_i)]]
={1}/{n} \sum_{i=1}^n[\|t\|^2-2\E[t(X_i)]]=\|t\|^2-2\int t f
=\|t-f\|^2-\|f\|^2$
and then minimizing $\gamma_n(t)$ comes down to minimizing the distance between $t$ and $f$.
So we define
\begin{equation}
  \label{estf}
  \hat{f}_m=\arg\min_{t\in S_m}\gamma_n(t)
\end{equation}
or, equivalently, 
$$\hat{f}_m=\sum_{j\in\mathbb{Z}}\hat{a}_j\varphi_{m,j}\quad\text{ with }
\quad \hat{a}_j=\frac{1}{n} \sum_{i=1}^n\u_{\varphi_{m,j}}(Y_i).$$
Actually we should define $\hat{f}_m=\sum_{|j|\leq K_n}\hat{a}_j\varphi_{m,j}$
because we can estimate only a finite number of coefficients. If $K_n$ is 
suitably chosen, it does not change the rate of convergence since the additional
terms can be made negligible. For the sake of simplicity, we let the sum over 
$\mathbb{Z}$. For an example of detailed truncation see \cite{comterota2006}.

Conditionally to $(X_i)$, the variance or stochastic error is 
$$\E[\|\hat{f}_m-f_m\|^2|X_1,\dots, X_n]\leq \sum_j\var[\frac{1}{n}\sum _{i=1}^n\u_{\varphi_{m,j}}(Y_i)|X_1,\dots, X_n]
\leq \frac{\|\sum_j\u_{\varphi_{m,j}}^2\|_\infty}{n}$$
since $Y_1,\dots, Y_n$ are independent conditionally to $(X_i)$. Then, it follows from Lemma 
\ref{ordrevar} (see Subsection \ref{techlem}) that $\|\sum_j\u_{\varphi_{m,j}}^2\|_\infty= \Delta(m)$ where
\begin{equation}
   \Delta(m)=\frac{1}{2\pi}\int_{-\pi m}^{\pi m}|q^*(u)|^{-2}du.
 \label{Deltam}\end{equation}
This implies that the order of the variance is $\Delta(m)/n.$
That is why we introduce $$\mathcal{M}_n=\left\{m\geq 1,\quad  \frac{\Delta(m)}{n}\leq 1\right\}.$$ 
To complete the estimation, we choose the best estimator among
the collection $(\hat{f}_m)_{m\in\mathcal{M}_n}$. Let 
$$\hat{m}=\underset{m\in\mathcal{M}_n}{\arg\min}\{\gamma_n(\hat{f}_m)+\pen(m)\}$$
where $\pen$ is  a penalty term to be specified later (see Theorem \ref{mainf}).
Finally we define $\tilde{f}=\hat{f}_{\hat{m}}$.

\subsection{Estimation of the density $F$ of $(X_i, X_{i+1})$}\label{estimationF}

We proceed similarly to the estimation of $f$. To define the contrast to minimize, we use the following lemma:
\begin{lem}\label{vt} For all function $T$, let $\v_T$ be the inverse Fourier transform of 
$T^*/(q^*\otimes q^*)(-.)$, i.e. $$\v_T(x,y)=
\frac{1}{4\pi^2}\iint e^{ixu+iyv} \frac{T^*(u,v)}{q^*(-u)q^*(-v)}dudv.$$ 
Then, for all $1\leq k\leq n$, 
  \begin{enumerate}
  \item $\E[\v_T(Y_k,Y_{k+1})|X_1,...,X_n]=T(X_k,X_{k+1})$
  \item $\E[\v_T(Y_k,Y_{k+1})]=\E[T(X_k,X_{k+1})]$
  \end{enumerate}
\end{lem}
For any function $T$ in $L^2(\mathbb{R}^2)$, we define the contrast 
$$\Gamma_n(T)=\frac{1}{n}\sum_{i=1}^n[\|T\|^2-2\v_T(X_i,X_{i+1})]$$
whose expectation is equal to $\|T\|^2-2/n\sum_{k=1}^n\E[T(X_k,X_{k+1})]
=\|T-F\|^2-\|F\|^2$.
We can now define
\begin{equation}
  \label{estF}
  \hat{F}_m=\arg\min_{T\in \mathbb{S}_m} \Gamma_n(T)
\end{equation}
By differentiating $\Gamma_n$, we obtain 
$$\hat{F}_m(x,y)=\sum_{j,k}\hat{A}_{j,k}\varphi_{m,j}(x)\varphi_{m,k}(y)
\text{ with } \hat{A}_{j,k}=\frac{1}{n}\sum_{i=1}^n\v_{\varphi_{m,j}\otimes \varphi_{m,k}}(Y_i, Y_{i+1}).$$
We choose again not to truncate the estimator for the sake of simplicity.
We have defined a collection of estimators $\{\hat{F}_m\}_{m\in \mathbb{M}_n}$. 
Note that  $\v_{t\otimes s}(x,y)=\u_t(x)\u_s(y)$ so that here
the variance is of order $\Delta^2(m)/n$, so we introduce
$$\mathbb{M}_n=\left\{m\geq 1,\quad  \frac{\Delta^2(m)}{n}\leq 1\right\}.$$
To define an adaptive estimator we have to select the best model $m$. So let
$$\hat{M}=\underset{m\in\mathbb{M}_n}{\arg\min}\{\Gamma_n(\hat{F}_m)+\PEN(m)\}$$
where $\PEN$ is a penalty function which is specified in Theorem \ref{mainF}.
Finally we consider the estimator $\tilde{F}=\hat{F}_{\hat{M}}$.

\subsection{Estimation of $\Pi$}\label{estimationpi}

Whereas the estimation of $f$ and $F$ is valid on the whole real line $\mathbb{R}$
or $\mathbb{R}^2$, we estimate $\Pi$ on a compact set $B^2$ only, because we need 
a lower bound on the stationary density. More precisely, we need to set some additional assumptions:
\begin{description}
  \item[A5] There exists a positive real $f_0$ such that $\forall x \in B, \quad f(x)\geq f_0$
  \item[A6] $\forall x \in B, \quad \forall y \in B,\qquad \Pi(x,y)\leq \|\Pi\|_{B,\infty}<\infty$
\end{description}
Now we set
\begin{equation}\label{defestpi}
  \tilde{\Pi}(x,y)=
\begin{cases}
  \cfrac{\tilde{F}(x,y)}{\tilde{f}(x)}& \text{ if } |\tilde{F}(x,y)|\leq n|\tilde{f}(x)|,\\
  0 &\text{ otherwise. }
\end{cases}\end{equation}
Here the truncation allows to avoid the too small values of $\tilde{f}$ in the quotient.
Now we evaluate upper bounds for the risk of our estimators.

\section{Results}

Our first theorem regards the problem of deconvolution. This result 
may be put together with results of \cite{comterota2006} in the i.i.d. case and of \cite{comtededeckertaupin2006}
in various mixing frameworks.

\begin{thm} \label{mainf}
Under Assumptions A1--A4, consider the estimator  $\tilde{f}=\hat{f}_{\hat{m}}$
where for each $m$,  $\hat{f}_m$ is defined by \eqref{estf} and 
$\hat{m}=\underset{m\in\mathcal{M}_n}{\arg\min}\{\gamma_n(\hat{f}_m)+\pen(m)\}$ with 
$$\pen(m)=k\cfrac{(\pi {m})^{[s-(1-s)_+/2]_+}\Delta(m)}{n}$$
where $k$ is a constant depending only on $k_0,k_1,b, \gamma, s$. Then there exists $C>0$ such that
$$\mathbb{E}\|\tilde{f}-f\|^2\leq 4\underset{m\in\mathcal{M}_n}{\inf}
\{\|f_m-f\|^2+ \pen(m)\}+\frac{C}{n}.$$
\end{thm}

The penalty is close to the variance order. It implies that the obtained rates of convergence are 
minimax in most cases. More precisely, the rates are given in the following corollary where 
$\lceil x\rceil$ denotes the ceiling function, i.e. the smallest integer larger than or equal to $x$.

\begin{cor}\label{ratesf}
  Under Assumptions of Theorem \ref{mainf},  if $f$ belongs to $\mathcal{A}_{\delta,r,a}(l)$, then
  \begin{itemize}
  \item If $r=0$ and $s=0$ $\qquad\displaystyle\mathbb{E}\|\tilde{f}-f\|^2\leq C n^{-\frac{2\delta}{2\delta+2\gamma+1}}$
  \item If $r=0$ and $s>0$ $\qquad\displaystyle\mathbb{E}\|\tilde{f}-f\|^2\leq C (\ln n)^{-2\delta/s}$
  \item If $r>0$ and $s=0$ $\qquad\displaystyle\mathbb{E}\|\tilde{f}-f\|^2\leq C \frac{(\ln n)^{(2\gamma+1)/r}}{n}$
  \item If $r>0$ and $s>0$
    \begin{itemize}
    \item if $r<s$ and  $k=\lceil (s/r-1)^{-1}\rceil - 1$
, there exist reals $b_i$
such that  $$\mathbb{E}\|\tilde{f}-f\|^2\leq C (\ln n)^{-2\delta/s}\exp[\sum_{i=0}^k b_i(\ln n)^{(i+1)r/s-i}]$$
    \item if $r=s$,   if $\xi=[2\delta b+(s-2\gamma-1-[s-(1-s)_+/2]_+)a]/[(a+b)s]$
$$\mathbb{E}\|\tilde{f}-f\|^2\leq Cn^{-{a}/(a+b)}(\ln n)^{-\xi}$$
    \item if $r>s$ and $k=\lceil (r/s-1)^{-1}\rceil - 1$, there exist reals $d_i$
such that $$\mathbb{E}\|\tilde{f}-f\|^2\leq C\frac{(\ln n)^{({1+2\gamma-s+[s-(1-s)_+/2]_+})/{r}}}
{n}\exp[-\sum_{i=0}^k d_i(\ln n)^{(i+1)s/r-i}]$$
    \end{itemize}
  \end{itemize}
\end{cor}

These rates are the same as those obtained in the case of i.i.d. variables $X_i$; they are studied 
in detail in \cite{comterota2006}. In the case $r>0, s>0$, we find the original rates obtained in 
\cite{lacourdeconvo}, proved as being optimal for $0<r<s$ in \cite{butuceatsybakov04}.
In the other cases, we can compare the 
 results of Theorem \ref{mainf} to the one obtained with a nonadaptive estimator. There is a loss 
 only in the case $r\geq s > 1/3$ where a logarithmic term is added. But in this case, the rates 
 are faster than any power of logarithm. 

Now let us study the risk for our estimator of the joined density $F$.

\begin{thm} \label{mainF}
Under Assumptions A1--A4, consider the estimator  $\tilde{F}=\hat{F}_{\hat{M}}$
where for each $m$,  $\hat{F}_m$ is defined by \eqref{estF} and 
$\hat{M}=\underset{m\in\mathbb{M}_n}{\arg\min}\{\Gamma_n(\hat{F}_m)+\PEN(m)\}$ with 
$$\PEN(m)=K\cfrac{(\pi {m})^{[s-(1-s)_+]_+}\Delta^2(m)}{n}$$
where $K$ is a constant depending only on $k_0,k_1,b, \gamma, s$. Then there exists $C>0$ such that
$$\mathbb{E}\|\tilde{F}-F\|^2\leq 4\underset{m\in\mathbb{M}_n}{\inf}
\{\|F_m-F\|^2+ \PEN(m)\}+\frac{C}{n}.$$
\end{thm}

The bases derived from the sine cardinal function are adapted to the estimation on the whole real line. 
The proof of Theorem \ref{mainF} actually contains the proof of another result 
(see Proposition~\ref{talagrand2} in Section \ref{proofs}):  the estimation 
of a bivariate density in a mixing framework on $\mathbb{R}^2$ and not only on a compact set. 
In this case of the absence of noise ($\varepsilon=0$), we obtain the same result with the penalty
$\PEN(m)=K_0(\sum_k\beta_{2k})m^2/n$. This limit case gives the mixing coefficients back in the penalty, 
as it always appears in this kind of estimation (see e.g. \cite{tribouleyviennet}).

It is then significant that in the presence of noise the penalty contains neither mixing term 
nor unknown quantity. It is entirely computable
since it depends only on the characteristic function $q^*$ of the noise which is known.

Theorem \ref{mainF} enables us to give rates of convergence for the estimation of $F$.

\begin{cor}\label{ratesF}
Under Assumptions of Theorem \ref{mainF}, if $F$ belongs to $\mathbb{A}_{\Delta,R,A}(L)$, then
  \begin{itemize}
  \item If $R=0$ and $s=0$ $\qquad\displaystyle\mathbb{E}\|\tilde{F}-F\|^2\leq C n^{-\frac{2\Delta}{2\Delta+4\gamma+2}}$
  \item If $R=0$ and $s>0$ $\qquad\displaystyle\mathbb{E}\|\tilde{F}-F\|^2\leq C (\ln n)^{-2\Delta/s}$
  \item If $R>0$ and $s=0$ $\qquad\displaystyle\mathbb{E}\|\tilde{F}-F\|^2\leq C \frac{(\ln n)^{(4\gamma+2)/R}}{n}$
  \item If $R>0$ and $s>0$
    \begin{itemize}
    \item if $R<s$ and  $k=\lceil (s/R-1)^{-1}\rceil - 1$, there exist reals $b_i$
such that  $$\mathbb{E}\|\tilde{F}-F\|^2\leq C (\ln n)^{-2\Delta/s}\exp[\sum_{i=0}^k b_i(\ln n)^{(i+1)R/s-i}]$$
    \item if $R=s$   if $\xi=[4\Delta b+(2s-4\gamma-2-[s-(1-s)_+]_+)A]/[(A+2b)s]$
$$\mathbb{E}\|\tilde{f}-f\|^2\leq Cn^{-{A}/(A+2b)}(\ln n)^{-\xi}$$
    \item if $R>s$ and $k=\lceil (R/s-1)^{-1}\rceil - 1$, there exist reals $d_i$
such that $$\mathbb{E}\|\tilde{F}-F\|^2\leq C\frac{(\ln n)^{({2+4\gamma-2s+[s-(1-s)_+]_+})/{R}}}
{n}\exp[-\sum_{i=0}^k d_i(\ln n)^{(i+1)s/R-i}]$$
    \end{itemize}
  \end{itemize}
\end{cor}

The rates of convergence look like the one of Corollary \ref{ratesf} with modifications
due to the bivariate nature of $F$. 
We can compare this result to the one of \cite{clemhmm} who studies only the case $R=0$ and $s=0$. He shows that the 
minimax lower bound in that case is $n^{-\frac{2\Delta}{2\Delta+4\gamma+2}}$, so our procedure is optimal,
whereas his estimator has a logarithmic loss for the upper bound.
We remark that if $s>0$ (supersmooth noise), the rate is logarithmic for $F$ belonging to a classic ordinary
smooth space. But if $F$ is also supersmooth, better rates are recovered.

Except in the case $R=0$ and $s=0$
, there is, to our knowledge, no lower bound available for this estimation.
We can however evaluate the performance of this estimator by comparing it with a nonadaptive estimator.
If the smoothness of $F$ is known, an $m$ depending on $R$ and $\Delta$ which minimizes 
the risk $\|F-F_m\|^2+\Delta(m)^2/n$ can be exhibited and then some rates of convergence for this nonadaptive estimator
are obtained.
As soon as  $s\leq 1/2$ (i.e. $[s-(1-s)_+]_+=0$), the penalty is $\Delta(m)^2/n$ and then 
the adaptive estimator recovers the same rates of convergence as those of a nonadaptive estimator
if the regularity of $F$ were known. It automatically minimizes the risk without prior knowledge on the regularity
of $F$ and there is no loss in the rates.
If $s>1/2$ a loss can appear but is not systematic. If $R<s$, the rate of convergence is unchanged since the bias
dominates. It is only in the case $R\geq s> 1/2$ that an additional logarithmic term appears. But in this case 
the risk decreases faster than any logarithmic function so that the loss is negligible.

We can now state the main result regarding the estimation of the transition density $\Pi$.

\begin{thm} \label{mainpi}
Under Assumptions A1--A6, consider the estimator  $\tilde{\Pi}$ 
defined in \eqref{defestpi}. We assume that $f$ belongs to $\mathcal{A}_{\delta,r,a}(l)$ 
with $\delta>1/2$ and 
that we browse only the models $m\in\mathcal{M}_n$ such that
\begin{equation}
  \label{eq}
  m \geq \ln\ln n\quad\text{ and }\quad m\Delta(m)\leq \frac{n}{(\ln n)^2}
\end{equation}
to define $\tilde{f}$.
 Then $\tilde{\Pi}$ verifies, for $n$ large enough,
$$\mathbb{E}\|\tilde{\Pi}-\Pi\|_B^2\leq C_1 \E\|\tilde{F}-F\|^2+C_2\E\|\tilde{f}-f\|^2
+\frac{C}{n}$$
where $\|T\|_B^2=\iint_{B^2} T^2(x,y)dxdy$.
\end{thm}

Note that, contrary to Theorems \ref{mainf} and \ref{mainF}, this result is asymptotic.
It states that the rate of convergence for $\Pi$ is no larger than the maximum of the rates of 
$f$ and $F$. The restrictions \eqref{eq} do not modify the conclusion of Theorem \ref{mainf}
and the resulting rates of convergence. Thus if $f$ and $F$ have the same regularity, the rates of convergence
for $\Pi$ are those of $F$, given in Corollary \ref{ratesF}.

If $s=0$ i.e. if $\varepsilon_i$ is ordinary smooth, then the rates of convergence 
are polynomial and even near the parametric rate $1/n$
if $R$ and $r$ are positive. But the smoother the error distribution is, the harder the estimation 
is. In the case of a supersmooth noise, the rates are logarithmic if $f$ or $F$ is ordinary smooth 
but faster than any power of logarithm if the hidden chain has supersmooth densities.
The exact rates depend on all regularities $\gamma$, $s$, 
$\delta$, $r$, $\Delta$, $R$ and are very tedious to write.
That is why we prefer to give some detailed examples.

\section{Examples}

\subsection{Autoregressive process of order 1}\label{ar1}

Let us study the case where the Markov chain is defined by 
$$X_{n+1}=\alpha X_n+\beta+\eta_{n+1}$$
where the $\eta_n$'s are i.i.d. centered Gaussian with variance $\sigma^2$.
This chain is irreducible, Harris recurrent and geometrically $\beta$-mixing.
The stationary distribution is Gaussian with mean $\beta/(1-\alpha)$ and 
variance $\sigma^2/(1-\alpha^2)$. So 
$$f^*(u)=\exp\left[-iu\left(\frac{\beta}{1-\alpha}\right)-\frac{\sigma^2}{2(1-\alpha^2)}u^2\right]$$
and then bias computing gives $\delta=1/2$, $r=2$. The function
$F$ is the density of a Gaussian vector with mean $(\beta/(1-\alpha),\beta/(1-\alpha))$ and 
variance matrix $\sigma^2/(1-\alpha^2)
  \begin{pmatrix}
1 & \alpha \\
\alpha & 1
  \end{pmatrix} $.
So 
$$F^*(u,v)=\exp\left[-i(u+v)\left(\frac{\beta}{1-\alpha}\right)-
\frac{\sigma^2}{2(1-\alpha^2)}(u^2+v^2+2\alpha uv)\right]$$
and  $\Delta=1/2$, $R=2$.

We can compute the rates of convergence for different kinds of noise $\varepsilon$.
If $\varepsilon$ has a Laplace distribution, $q^*(u)=1/(1+u^2)$ so $s=0$, $\gamma=2$.
In this case, Corollary \ref{ratesf} gives $\mathbb{E}\|\tilde{f}-f\|^2\leq C (\ln n)^{5/2}/{n}$
and $\mathbb{E}\|\tilde{F}-F\|^2\leq C {(\ln n)^5}/{n}$. Consequently,
$$\mathbb{E}\|\tilde{\Pi}-\Pi\|_B^2\leq C \frac{(\ln n)^5}{n}$$
with $B$ an interval $[-d,d]$.
This rate is close to the parametric rate $1/n$; it is due to the great smoothness of the chain 
compared with that of error.


If now $\varepsilon$ has a normal distribution with variance $\tau^2$, then we compute
$$\mathbb{E}\|\tilde{\Pi}-\Pi\|_B^2\leq C n^{-\frac{\sigma^2}{\sigma^2+2\tau^2}}(\ln n)^{-\frac{\tau^2}{\sigma^2+2\tau^2}}.$$

\subsection{Cox-Ingersoll-Ross process}

Another example is given by $X_n=R_{n\tau}$ with $\tau$ a fixed sampling interval and 
$R_t$ the so-called Cox-Ingersoll-Ross process defined by 
$$dR_t=(2\theta R_t+\kappa\sigma_0^2)dt +2 \sigma_0\sqrt{R_t}dW_t\qquad\qquad \theta<0, \kappa\in\{2,3,...\}.$$
Following \cite{chaleyatgenon06}, we observe that $X_n$ is the square of the Euclidean norm of 
a $\kappa$-dimensional vector whose components are linear autoregressive processes of order 1.
The stationary distribution is a Gamma distribution with parameter $\kappa/2$ and $|\theta|/\sigma^2$ so that 
$$f^*(u)=\left(1+iu\frac{\sigma_0^2}{|\theta|}\right)^{-\kappa/2}$$
and $r=0$, $\delta=(\kappa-1)/2$.
To compute the characteristic function of the joined density, we write
$$F^*(u,v)=\int \E[e^{-ivX_1}|X_0=x]e^{-iux}f(x)dx.$$
Let $\beta^2=\sigma_0^2(e^{2\theta\tau}-1)/(2\theta)$. Then, 
conditionally to $X_0=x$, $\beta^{-2}X_1$ is a non-central chi-square
$\chi'^2(e^{2\theta\tau}x/\beta^2,\kappa)$, so that
$$\E[e^{-ivX_1}|X_0=x]=(1+2iv\beta^2)^{-\kappa/2}
\exp\left(-\cfrac{ive^{2\theta\tau}x}{1+2iv\beta^2}\right).$$
This implies
$$F^*(u,v)=\left[1-(1-e^{2\theta\tau})\frac{\sigma_0^4}{\theta^2}uv+i\frac{\sigma_0^2}{|\theta|}(u+v)\right]^{-\kappa/2}$$
and $R=0$, $\Delta=(\kappa-1)/2$.
Then, if for example the noise has a Gaussian distribution ($\gamma=0, s=2$),
the rate of convergence is $(\ln n)^{(1-\kappa)/2}$. But this rate is faster if $\varepsilon$
has a Gamma distribution with shape parameter $\alpha$ (so that $\gamma=\alpha$, $s=0$): we obtain 
in this case $n^{(1-\kappa)/(\kappa +4\alpha+1)}$.

\subsection{Stochastic volatility model}

Our work allows to study some multiplicative models as the so-called stochastic volatility model in finance (see
\cite{genonjantheaularedo2000} for the links between the standard continuous-time SV models and the hidden Markov models). Let
$$Z_n=U_n^{1/2}\eta_n$$
where $(U_n)$ is a nonnegative Markov chain, $(\eta_n)$ a sequence of i.i.d. standard Gaussian variables, the two 
sequences being independent. Setting 
$X_n=\ln (U_n)$ and $\varepsilon_n=\ln (\eta_n^2)$ leads us back to our initial problem.

The noise distribution is the logarithm of a chi-square distribution
and then verifies 
$q^*(x)=2^{-ix}\Gamma(1/2-ix)/\sqrt{\pi}$. \cite{vanesspreij2005} show that
$|q^*(x)|\sim_{+\infty}\sqrt{2}e^{-\pi|x|/2}$ and then $s=1,\gamma=0$.


We assume that the logarithm of the hidden chain $X_n$ derives from a regular sampling of an 
Ornstein-Uhlenbeck process, i.e. $X_n=V_{n\tau}$ where $V_t$ is defined by the equation $$dV_t=\theta V_tdt +\sigma dB_t$$
with $B_t$ a standard Brownian motion. Then all the assumptions are satisfied. Similarly to Subsection \ref{ar1},
the stationary distribution is Gaussian with mean $0$ and variance $\sigma^2/2|\theta|$
and then $\delta=1/2$, $r=2$. In the same way $F$ is the density of a centered Gaussian vector with variance matrix
$\sigma^2/(2|\theta|)
  \begin{pmatrix}
1 & e^{ \theta\tau}\\
e^{ \theta\tau} & 1
  \end{pmatrix} $  and then $\Delta=1/2$, $R=2$.
We  obtain the following rate of convergence on some interval $B=[-d,d]$
$$\mathbb{E}\|\tilde{\Pi}-\Pi\|_B^2\leq C \sqrt{\ln n}\frac{\exp[(\pi/\beta)\sqrt{\ln n}]}{n}$$
with $\beta^2=\sigma^2(e^{2\theta\tau}-1)/(2\theta)$.

\section{Proofs}\label{proofs}


Here we do not prove the results concerning the estimation of $f$. Indeed they are similar to the
ones concerning $F$ (but actually simpler) and the ones of \cite{comterota2006}. It is then sufficient to use 
corresponding proofs for $F$ mutatis mutandis. 

For the sake of simplicity, all constants in the following are denoted by $C$, even
if they have different values.

\subsection{Proof of Lemma \ref{vt}}
It is sufficient to prove the first assertion.
First we write that
$\v_T(Y_k,Y_{k+1})={1}/{4\pi^2}\int e^{iY_ku+iY_{k+1}v} {T^*(u,v)}/{q^*(-u)q^*(-v)}dudv$
so that
$$\E[\v_T(Y_k,Y_{k+1})|X_1,...,X_n]=\frac{1}{4\pi^2}\int \E[e^{iY_ku+iY_{k+1}v}|X_1,...,X_n] 
\frac{T^*(u,v)}{q^*(-u)q^*(-v)}dudv.$$
By using the independence between $(X_i)$ and $(\varepsilon_i)$, we compute
\begin{eqnarray*}
&&  \E[e^{iY_ku+iY_{k+1}v}|X_1,..,X_n]=\E[e^{iX_ku+iX_{k+1}v}e^{i\varepsilon_ku+i\varepsilon_{k+1}v}
|X_1,..,X_n]\\ && =e^{iX_ku+iX_{k+1}v}\E[e^{i\varepsilon_ku}]\E[e^{i\varepsilon_{k+1}v}]
=e^{iX_ku+iX_{k+1}v}\int e^{ixu} q(x)dx\int e^{iyv} q(y)dy\\ &&
=e^{iX_ku+iX_{k+1}v}q^*(-u)q^*(-v).
\end{eqnarray*}
Then 
\begin{eqnarray*}
\E[\v_T(Y_k,Y_{k+1})|X_1,..,X_n]&=&\frac{1}{4\pi^2}\int e^{iX_ku+iX_{k+1}v}q^*(-u)q^*(-v)
\frac{T^*(u,v)}{q^*(-u)q^*(-v)}dudv\\
&=&\frac{1}{4\pi^2}\int e^{iX_ku+iX_{k+1}v}{T^*(u,v)}dudv=T(X_k,X_{k+1}).
\end{eqnarray*}

\subsection{Proof of Theorem \ref{mainF}}

First we introduce some auxiliary variables whose existence is ensured by Assumption A4  
of mixing. In the case of arithmetical mixing, since $\theta>6$, there exists a real $c$ such that 
$0<c<1/2$ and $c \theta>3$.
We set in this case $q_n=\frac{1}{2}\lfloor n^c \rfloor.$  
In the case of geometrical mixing, we set $q_n=\frac{1}{2}\lfloor c\ln(n) \rfloor$ 
where $c$ is a real larger than $3/\theta$.

For the sake of simplicity, we suppose that $n=4p_nq_n$, with $p_n$ an integer.
Let for $i=1,\dots,n/2$, $V_i=(X_{2i-1},X_{2i})$
 and for $l=0,\dots,p_n-1$, $A_l= (V_{2lq_n+1},..., V_{(2l+1)q_n})$, $B_l=(V_{(2l+1)q_n+1},..., V_{(2l+2)q_n})$.
As in \cite{viennet97}, by using Berbee's coupling Lemma, we can build a sequence $(A_l^*)$ such that 
$$\begin{cases}
A_l \text{ and } A_l^* \text{ have the same distribution,}\\
A_l^* \text { and } A_{l'}^* \text{ are independent if }l\neq l',\\
P(A_l\neq A_{l}^*)\leq \beta_{2q_n}.
\end{cases}$$
In the same way, we build $(B_l^*)$ and we define for any $l\in\{0,\dots, p_n-1\}$, \\
$A_l^*= (V^*_{2lq_n+1},..., V^*_{(2l+1)q_n})$, $B_l^*= (V^*_{(2l+1)q_n+1},..., V^*_{(2l+2)q_n})$
so that the sequence $(V_1^*,\dots, V_{n/2}^*)$ and then the sequence $(X_1^*,\dots, X_{n}^*)$
are well defined.
We can now define 
$$\Omega^*=\{\forall i, 1\leq i \leq n \quad X_i=X_i^*\}.$$
Then we split the risk into two terms:
\begin{equation*}
  \E(\|\tilde{F}-F\|^2)=\E(\|\tilde{F}-F\|^2\1_{\Omega^{*}})+\E(\|\tilde{F}-F\|^2\1_{\Omega^{*c}}).
\end{equation*}
To pursue the proof, we observe that for all $T,T'$
$$\Gamma_n(T)-\Gamma_n(T')=\|T-F\|^2-\|T'-F\|^2 -2 Z_n(T-T')$$
where $\displaystyle Z_n(T)=\frac{1}{n}\sum_{i=1}^{n}\left\{
\v_T(Y_i,Y_{i+1})-\int T(x,y) F(x,y)dxdy\right\}.$\\
Let us fix $m\in\mathbb{M}_n$ and denote by ${F}_m$ the orthogonal projection of $F$ on $\mathbb{S}_m$.
Since $\Gamma_n(\tilde{F})+\PEN(\hat{M})\leq \Gamma_n({F}_m)+\PEN(m)$, we have
\begin{eqnarray*}
\|\tilde{F}-F\|^2&\leq& \|F_m-F\|^2+2Z_n(\tilde{F}-F_m)+\PEN(m)-\PEN(\hat{M})\\
&\leq&\|F_m-F\|^2  +2\|\tilde{F}-F_m\|\sup_{T\in B(\hat{M})}Z_n(T)+\PEN(m)-\PEN(\hat{M})
\end{eqnarray*}
where, for all $m'$, $B(m')=\{T\in \S_m+\S_{m'},\quad \|T\|=1\}$. Then, 
using inequality $2xy\leq x^2/4+4y^2$,
\begin{equation}\label{ab}
\|\tilde{F}-F\|^2\leq \|F_m-F\|^2+\frac{1}{4}\|\tilde{F}-F_m\|^2+4\sup_{T\in B(\hat{M})}Z_n^2(T)
+\PEN(m)-\PEN(\hat{M}).
\end{equation}
By denoting $\E_X$ the expectation conditionally to $X_1,\dots, X_n$ and by using
Lemma~\ref{vt},  $Z_n(T)$ can be split into two terms :
$$Z_n(T)=Z_{n,1}(T)+Z_{n,2}(T)$$ with 
$$\begin{cases} 
\displaystyle Z_{n,1}(T)=\frac{1}{n}\sum_{i=1}^{n}\left\{
\v_T(Y_i,Y_{i+1})-\E_X[\v_T(Y_i,Y_{i+1})]\right\},\\
\displaystyle Z_{n,2}(T)=\frac{1}{n}\sum_{i=1}^{n}\left\{
T(X_i,X_{i+1})-\int T(x,y) F(x,y)dxdy\right\}.\\
\end{cases}$$
Now let $P_1(.,.)$ be a function such that for all $m,m'$,
\begin{equation}
  \label{pmm'}
  16 P_1(m,m')\leq \PEN(m)+\PEN(m').
\end{equation}
 Then \eqref{ab} becomes
\begin{eqnarray*}
\|\tilde{F}-F\|^2\leq \|F_m-F\|^2+\frac{1}{2}(\|\tilde{F}-F\|^2+\|F-F_m\|^2)+2\PEN(m)\\ 
+8[\sup_{T\in B(\hat{M})}Z_{n,1}^2(T)-P_1(m,\hat{M})]+8[\sup_{T\in B(\hat{M})}Z_{n,2}^2(T)-P_1(m,\hat{M})]
\end{eqnarray*}
which gives, by introducing a function $P_2(.,.)$,
\begin{eqnarray*}
\frac{1}{2}\|\tilde{F}-F\|^2\1_{\Omega^*}\leq \frac{3}{2}\|F_m-F\|^2+2\PEN(m)
+8\sum_{m'\in\mathbb{M}_n}[\sup_{T\in B(m')}Z_{n,1}^2(T)-P_1(m,m')]_+\\
+8\sum_{m'\in\mathbb{M}_n}[\sup_{T\in B(m')}Z_{n,2}^2(T)-P_2(m,m')]_+\1_{\Omega^*}
+8\sum_{m'\in\mathbb{M}_n}[P_2(m,m')-P_1(m,m')].
\end{eqnarray*}

We now use the following propositions:

\begin{prop}\label{talagrand1}
Let $P_1(m,m')=C(q)(\pi {m''})^{[s-(1-s)_+]_+}\Delta^2(m'')/n$ where $\Delta(m)$ is defined in 
\eqref{Deltam} and $m''=\max(m,m')$ and C(q) is a constant. Then, under assumptions of Theorem 
\ref{mainF}, there exists a positive constant $C$ such that  \begin{equation}\sum_{m'\in \mathbb{M}_n}
\E\left(\left[\sup_{T\in B(m')}Z_{n,1}^2(T)-P_1(m,m')\right]_+\right)\leq \frac{C}{n}.
\label{inegtalagrand1}\end{equation}
\end{prop}

\begin{prop}\label{talagrand2}
Let $P_2(m,m')=96(\sum_k\beta_{2k}){{m''}}/{n}$ where $\Delta(m)$ is defined in 
\eqref{Deltam} and $m''=\max(m,m')$. Then, under assumptions of Theorem \ref{mainF},
there exists a positive constant $C$ such that  \begin{equation}\sum_{m'\in \mathbb{M}_n}
\E\left(\left[\sup_{T\in B(m')}Z_{n,2}^2(T)-P_2(m,m')\right]_+\1_{\Omega^*}\right)\leq \frac{C}{n}.
\label{inegtalagrand2}\end{equation}
\end{prop}

The definition of the functions $P_1(m,m')$ and $P_2(m,m')$ given in Propositions 
\ref{talagrand1} and~\ref{talagrand2} imply that there exists $m_0$ 
such that $\forall m'> m_0 \quad P_1(m,m')\geq P_2(m,m')$. 
(If $s=0=\gamma$ (case of a null noise), it would be wrong and the penalty would then be 
$P_2(m,m')$ instead of $ P_1(m,m'$)).
Then 
\begin{eqnarray}\label{cc}
\sum_{m'\in\mathbb{M}_n}[P_2(m,m')-P_1(m,m')]\leq\sum_{m'\leq  m_0}P_2(m,m')
\leq \frac{C(m_0)}{n}.
\end{eqnarray}

Since ${m''}\Delta^2(m'')\leq {m}\Delta^2(m)+{m'}\Delta^2(m')$, the condition 
\eqref{pmm'} is verified with 
$$\PEN(m)=16C(q)(\pi {m})^{[s-(1-s)_+]_+}\frac{\Delta^2(m)}{n}
.$$And finally, combining \eqref{cc} and Propositions \ref{talagrand1} and~\ref{talagrand2},
\begin{eqnarray*}
\E(\|\tilde{F}-F\|^2\1_{\Omega^*})\leq 4 (\|F_m-F\|^2+\PEN(m))+\frac{C}{n}.
\end{eqnarray*}

For the term $\E(\|\tilde{F}-F\|^2\1_{\Omega^{*c}})$, recall that 
$$\hat{F}_m(x,y)=\sum_{j,k}\hat{A}_{j,k}\varphi_{m,j}(x)\varphi_{m,k}(y)
\text{ with } \hat{A}_{j,k}=\frac{1}{n}\sum_{i=1}^n\v_{\varphi_{m,j}\otimes \varphi_{m,k}}(Y_i, Y_{i+1}).$$
Thus, for any $m$ in $\mathbb{M}_n$, \begin{eqnarray}
\|\hat{F}_m\|^2&=&\sum_{j,k}\left[\frac{1}{n}\sum_{i=1}^n\v_{\varphi_{m,j}\otimes \varphi_{m,k}}(Y_i, Y_{i+1})\right]^2
\leq \frac{1}{n^2}\sum_{j,k}n\sum_{i=1}^n\v_{\varphi_{m,j}\otimes \varphi_{m,k}}^2(Y_i, Y_{i+1})\nonumber\\
&\leq& \|\sum_{j,k}\v_{\varphi_{m,j}\otimes \varphi_{m,k}}^2\|_\infty\leq 
\|\sum_{j}\u_{\varphi_{m,j}}^2\|_\infty^2\leq \Delta^2(m)\label{dd}
\end{eqnarray}
using Lemma \ref{ordrevar}.
Then 
$\|\hat{F}_{\hat{M}}\|^2\leq \Delta^2(\hat{M})\leq n$ since $\hat{M}$ belongs 
to $\mathbb{M}_n$. And $$\E\|\tilde{F}-F\|^2\1_{\Omega^{*c}}\leq 
\E(2(\|\tilde{F}\|^2+\|F\|^2)\1_\Omega^{*c})\leq 2(n+\|F\|^2)P(\Omega^{*c}).$$
Using Assumption A4 in the geometric case, $\beta_{2q_n}\leq Me^{-\theta c\ln(n)}\leq M n^{-\theta c}$ 
and, in the other case, $\beta_{2q_n}\leq M(2q_n)^{-\theta}\leq M n^{-\theta c}$. 
Then $P(\Omega^{*c})\leq 2p_n\beta_{2q_n}\leq nMn^{-c\theta}.$
Since $c\theta>3$, $P(\Omega^{*c})\leq{M}n^{-2}$, which implies
$E(\|\tilde{F}-F\|^2\1_{\Omega^{*c}})\leq {C}/{n}$.

Finally we obtain \begin{eqnarray*}
\E\|\tilde{F}-F\|^2&\leq &\E(\|\tilde{F}-F\|^2\1_{\Omega^*})+E(\|\tilde{F}-F\|^2\1_{\Omega^{*c}})\\
&\leq &4 (\|F_m-F\|^2+\PEN(m))+\frac{C}{n}.
\end{eqnarray*}
This inequality holds for each $m\in \mathbb{M}_n$, so the result is proved.

\subsection{Proof of Proposition \ref{talagrand1}}

We start by isolating odd terms from even terms to avoid overlaps:
$$Z_{n,1}(T)=\frac{1}{2}Z_{n,1}^o(T)+\frac{1}{2}Z_{n,1}^e(T)$$ with 
$$\begin{cases} 
 \displaystyle Z_{n,1}^o(T)=\frac{2}{n}\sum_{i=1, i \text{ odd}}^{n}\left\{
\v_T(Y_i,Y_{i+1})-\E_X[\v_T(Y_i,Y_{i+1})]\right\},\\
\displaystyle Z_{n,1}^e(T)=\frac{2}{n}\sum_{i=1, \text{ even }}^{n}\left\{
\v_T(Y_i,Y_{i+1})-\E_X[\v_T(Y_i,Y_{i+1})]\right\}.\\
\end{cases}$$
It is sufficient to deal with the first term only, as the second one is similar.
For each $i$, let $U_i=(Y_{2i-1},Y_{2i})$, then
$$Z_{n,1}^o(T)=\frac{1}{n/2}\sum_{i=1}^{n/2}\left\{\v_T(U_i)-\E_X[\v_T(U_i)]\right\}.$$
Notice that conditionally to $X_1,\dots, X_n$, the $U_i$'s are independent. Thus we can use the 
Talagrand inequality recalled in Lemma \ref{tala}. 
Note that if $T$ belongs to $\mathbb{S}_m+\mathbb{S}_{m'}$, then $T$ can be written
$T_1+T_2$ where $T_1^*$ has its support in $[-\pi {m},\pi {m}]^2$
and $T_2^*$ has its support in $[-\pi {m'},\pi {m'}]^2$. 
Then  $T$ belongs to $\mathbb{S}_{m''}$ where $m''$ is defined by 
\begin{equation}\label{m''}
{m''}=\max (m,{m'}).
\end{equation} Now let us compute $M_1$, $H$ and $v$ of the Talagrand's inequality.
\begin{enumerate}

\item If $T$ belongs to $B(m')$, $$\v_T(x,y)=\sum_{j,k}a_{jk}\v_{\varphi_{m'',j}\otimes \varphi_{m'',k}}(x,y)
=\sum_{j,k}a_{jk}\u_{\varphi_{m'',j}}(x)\u_{\varphi_{m'',k}}(y).$$ Thus 
$|\v_T(x,y)|^2\leq \sum_{j,k}|\u_{\varphi_{m'',j}}(x)\u_{\varphi_{m'',k}}(y)|^2$.
So $$\sup_{T\in B(m')}\|\v_T\|_\infty^2\leq \|\sum_{j,k}|\u_{\varphi_{m'',j}}(x)\u_{\varphi_{m'',k}}(y)|^2\|_\infty
\leq \|\sum_{j}|\u_{\varphi_{m'',j}}|^2\|_\infty^2.$$
By using Lemma \ref{ordrevar}, $M_1=\Delta(m'').$

\item To compute $H^2$, we write 
  \begin{eqnarray*}
 \E_X(\sup_{T\in B(m')}(Z_{n,1}^{o})^{2}(T))&\leq &
\E_X(\sum_{j,k}Z_{n,1}^{o}(\varphi_{m'',j}\otimes \varphi_{m'',k})^2)  \\ 
&\leq& \sum_{j,k}\var_X\left[\frac{2}{n}\sum_{i=1, i \text{ odd}}^n
\u_{\varphi_{m'',j}}(Y_i)\u_{\varphi_{m'',k}}(Y_{i+1})\right]\\
&\leq &\sum_{j,k}\frac{4}{n^2}\sum_{i=1,i \text{ odd}}^n
\var_X[\u_{\varphi_{m'',j}}(Y_i)\u_{\varphi_{m'',k}}(Y_{i+1})]
  \end{eqnarray*}
since, conditionally to $X_1,\dots, X_n$, the $U_i$'s are independent. And then
  \begin{eqnarray*}
\E_X(\sup_{T\in B(m')}Z_{n,1}^{o2}(T))\leq \sum_{j,k}\frac{4}{n^2}\sum_{i=1,i \text{ odd}}^n
\E_X[\u_{\varphi_{m'',j}}^2(Y_i)\u_{\varphi_{m'',k}}^2(Y_{i+1})]\\
\leq \frac{4}{n^2}\sum_{i=1,i \text{ odd}}^n\|\sum_{j}|\u_{\varphi_{m'',j}}|^2\|_\infty
\|\sum_{k}|\u_{\varphi_{m'',k}}|^2\|_\infty\leq \frac{2\Delta(m'')^2}{n}.
\end{eqnarray*}So we set $H=\sqrt{2}{\Delta(m'')}/{\sqrt{n}}.$

\item We still have to find $v$. On the one hand
  \begin{eqnarray*}
\var_X[\v_T(Y_k, Y_{k+1})]\leq \E_X[(\sum_{j,k}a_{jk}
\u_{\varphi_{m'',j}}(Y_k)\u_{\varphi_{m'',k}}(Y_{k+1}))^2]\\
\leq \sum_{j,k}a_{jk}^2\|\sum_{j}|\u_{\varphi_{m'',j}}|^2\|_\infty
\|\sum_{k}|\u_{\varphi_{m'',k}}|^2\|_\infty
\end{eqnarray*}
and so $v\geq \Delta(m'')^2$.
On the other hand
\begin{eqnarray}
 && \var_X[\v_T(Y_k, Y_{k+1})]\nonumber\\&&\leq \sum_{j_1,k_1}\sum_{j_2,k_2}a_{j_1k_1}a_{j_2k_2}
\E_X[\u_{\varphi_{m'',j_1}}\u_{\varphi_{m'',j_2}}(Y_k)\u_{\varphi_{m'',k_1}}
\u_{\varphi_{m'',k_2}}(Y_{k+1})]\nonumber\\&&
\leq \sum_{j,k}a_{jk}^2\sqrt{\sum_{j_1,k_1}\sum_{j_2,k_2}\E_X^2[\u_{\varphi_{m'',j_1}}
\u_{\varphi_{m'',j_2}}(Y_k)\u_{\varphi_{m'',k_1}}\u_{\varphi_{m'',k_2}}(Y_{k+1})]}\nonumber\\
&&\leq \sum_{j,k}a_{jk}^2\sqrt{\sum_{j_1,j_2}\E_X^2[\u_{\varphi_{m'',j_1}}
\u_{\varphi_{m'',j_2}}(Y_k)]\sum_{k_1,k_2}\E_X^2[\u_{\varphi_{m'',k_1}}
\u_{\varphi_{m'',k_2}}(Y_{k+1})]},\label{aa}
\end{eqnarray}
using conditional independence. Now we use Lemma \ref{ordrevar} to compute 
\begin{eqnarray*}
\E_X[\u_{\varphi_{m'',j_1}}\u_{\varphi_{m'',j_2}}(Y_k)]=\int (\u_{\varphi_{m'',j_1}}\u_{\varphi_{m'',j_2}})(X_k+x)q(x)dx\\
=\frac{{m''}}{4\pi^2}\int\int_{-\pi}^{\pi} \frac{e^{-ij_1v}e^{i(x+X_k)v{m''}}}{q^*(-v{m''})}dv
\int_{-\pi}^{\pi} \frac{e^{-ij_2u}e^{i(x+X_k)u{m''}}}{q^*(-u{m''})}duq(x)dx\\
=\frac{{m''}}{4\pi^2}\int_{-\pi}^{\pi} \int_{-\pi}^{\pi} \frac{e^{-ij_1v-ij_2u}
e^{iX_k(u+v){m''}}}{q^*(-v{m''})q^*(-u{m''})}\int e^{ix(u+v){m''}}q(x)dxdudv.
\end{eqnarray*}
If we set $W(u,v)={m''}{e^{iX_k(u+v){m''}}q^*(-(u+v){m''})}/[{q^*(-v{m''})q^*(-u{m''})}]$, then
$\E_X[\u_{\varphi_{m'',j_1}}\u_{\varphi_{m'',j_2}}(Y_k)]$ is the Fourier coefficient with order 
$(j_1,j_2)$ of $W$. Using Parseval's formula
\begin{eqnarray*}
  \sum_{j_1,j_2}\E_X^2[\u_{\varphi_{m'',j_1}}\u_{\varphi_{m'',j_2}}(Y_k)]=\frac{1}{4\pi^2}
\int_{-\pi}^{\pi}\int_{-\pi}^{\pi}|W(u,v)|^2dudv\\
=\frac{{m''}^2}{4\pi^2}\int_{-\pi}^{\pi}\int_{-\pi}^{\pi}\left|\frac{q^*(-(u+v){m''})}{q^*(-v{m''})q^*(-u{m''})}\right|^2dudv.
\end{eqnarray*}
Now we apply Schwarz inequality:
\begin{eqnarray*}
 && \sum_{j_1,j_2}\E_X^2[\u_{\varphi_{m'',j_1}}(Y_k)]\leq \\ &&\frac{{m''}^2}{4\pi^2}
\sqrt{\iint \frac{|q^*(-(u+v){m''})|^2}{|q^*(-u{m''})|^4}dudv}
\sqrt{\iint \frac{|q^*(-(u+v){m''})|^2}{|q^*(-v{m''})|^4}dudv}\\ &&
\leq \frac{{m''}}{4\pi^2}\int_{-\pi}^{\pi}|q^*(-u{m''})|^{-4}du \int |q^*(x)|^2dx
\leq \frac{\|q\|^2}{4\pi^2}\int_{-\pi {m''}}^{\pi {m''}}|q^*(-u)|^{-4}du .
\end{eqnarray*}
We introduce the following notation:
\begin{equation}
\Delta_2(m)=\frac{1}{4\pi^2}\int_{-\pi m}^{\pi m}|q^*(u)|^{-4}du.
  \label{Delta2m}
 \end{equation}
Finally, coming back to \eqref{aa},
$\var_X[\v_T(Y_k, Y_{k+1})]\leq \|T\|^2\|q\|^2\Delta_2(m'')$
which yields $v\geq \|q\|^2\Delta_2(m'')$.
Finally we write $v=\min( \|q\|^2\Delta_2(m''), \Delta^2(m'')).$
\end{enumerate}

We can now use Talagrand's inequality (see Lemma \ref{tala}):
\begin{eqnarray*}
  \E[\sup_{T\in B(m')}(Z_{n,1}^{o})^{2}(T)-2(1+2\epsilon)\frac{2\Delta^2(m'')}{n}]_+\leq \\
\frac{C}{n}\{ve^{-K_1\epsilon\Delta^2(m'')/v}+\frac{\Delta^2(m'')}{nC^2(\epsilon)}
e^{-K_2C(\epsilon)\sqrt{\epsilon}\sqrt{n}}\}.
\end{eqnarray*}
And then, if $P_1(m,m')\geq 4(1+2\epsilon){\Delta^2(m'')}/n$,
$$\sum_{m'\in\mathbb{M}_n}\E\left[\sup_{T\in B(m')}(Z_{n,1}^{o})^{2}(T)-P_1(m,m')\right]_+
\leq \frac{K}{n}\{I(m)+II(m)\}$$
with $I(m)=\sum_{m'\in\mathbb{M}_n}ve^{-K_1\epsilon\Delta^2(m'')/v}$; $II(m)=\sum_{m'\in\mathbb{M}_n}({1}/{nC^2(\epsilon)})
e^{-K_2C(\epsilon)\sqrt{\epsilon}\sqrt{n}}.$

To bound these terms, we use Lemma \ref{del} which yields to 
$$v\leq c_3 (\pi {m''})^{4\gamma+\min(1-s,2-2s)}e^{4b(\pi {m''})^s}\text{ and  }
\frac{\Delta^2({m''})}{v}\geq c_4 (\pi {m''})^{(1-s)_+}$$
where $c_3$ and $c_4$ depend only on $k_0,k_1,\gamma$ and $s$.
Therefore,
\begin{eqnarray*}
  I(m)&\leq &c_3\sum_{m'\in\mathbb{M}_n}(\pi {m''})^{4\gamma+\min(1-s,2-2s)}
e^{4b(\pi {m''})^s-K_1c_4\epsilon (\pi {m''})^{(1-s)_+}}\\
&\leq &c_3\sum_{m'\in\mathbb{M}_n}[(\pi {m})^{4\gamma+\min(1-s,2-2s)}
e^{4b(\pi {m})^s}\\ &&+(\pi {m'})^{4\gamma+\min(1-s,2-2s)}e^{4b(\pi {m'})^s}]
e^{-\frac{K_1c_4\epsilon }{2}[(\pi {m})^{(1-s)_+}+(\pi {m'})^{(1-s)_+}]}\\
&\leq &c_3(\pi {m})^{4\gamma+\min(1-s,2-2s)}e^{4b(\pi {m})^s
-\frac{K_1c_4\epsilon }{2}(\pi {m})^{(1-s)_+}}
\sum_{m'\in\mathbb{M}_n}e^{-\frac{K_1c_4\epsilon }{2}(\pi {m'})^{(1-s)_+}}\\ &&
+c_3e^{-\frac{K_1c_4\epsilon }{2}(\pi {m})^{(1-s)_+}}
\sum_{m'\in\mathbb{M}_n}(\pi {m'})^{4\gamma+\min(1-s,2-2s)}
e^{4b(\pi {m'})^s-\frac{K_1c_4\epsilon }{2}(\pi {m'})^{(1-s)_+}}.
\end{eqnarray*}
We have to distinguish three cases
\begin{description}

\item[case $s<(1-s)_+\Leftrightarrow s<1/2$]
In this case we choose $\epsilon=8b/(K_1c_4)$ and then 
\begin{eqnarray*}
  I(m)&\leq &c_3(\pi {m})^{4\gamma+1-s}e^{4b[(\pi {m})^s- (\pi {m})^{(1-s)}]}
\sum_{m'\in\mathbb{M}_n}e^{-K_1c_4(\pi {m'})^{(1-s)}}\\ &&
+c_3e^{-K_1c_4(\pi {m})^{(1-s)}}\sum_{m'\in\mathbb{M}_n}(\pi {m'})^{4\gamma+1-s}
e^{4b[(\pi {m'})^s-(\pi {m'})^{(1-s)}]}
\end{eqnarray*}
which implies that $I(m)$ is bounded. Moreover the definition of $\mathbb{M}_n$ and 
Lemma \ref{del} give $ |\mathbb{M}_n|\leq Cn^\zeta$ with $C>0$ and $\zeta>0$. So
$II(m)\leq (C/n)|\mathbb{M}_n|e^{-K'_2\sqrt{n}} $ is bounded too.

\item[case $s=(1-s)_+\Leftrightarrow s=1/2$]
In this case
\begin{eqnarray*}
  I(m)&\leq &c_3(\pi {m})^{4\gamma+1/2}e^{(4b-\frac{K_1c_4\epsilon }{2})(\pi {m})^{1/2}}
\sum_{m'\in\mathbb{M}_n}e^{-\frac{K_1c_4\epsilon }{2}(\pi {m'})^{1/2}}\\ &&
+c_3e^{-\frac{K_1c_4\epsilon }{2}(\pi {m})^{1/2}}
\sum_{m'\in\mathbb{M}_n}(\pi {m'})^{4\gamma+1/2}
e^{(4b-\frac{K_1c_4\epsilon }{2})(\pi {m'})^{1/2}}.
\end{eqnarray*}
We choose $\epsilon$ such that $4b-{K_1c_4\epsilon}/{2}=-4b$ so that
\begin{eqnarray*}
  I(m)&\leq &c_3(\pi {m})^{4\gamma+1/2}e^{-4b(\pi {m})^{1/2}}
\sum_{m'\in\mathbb{M}_n}e^{-\frac{K_1c_4\epsilon }{2}(\pi {m'})^{1/2}}\\ &&
+c_3e^{-\frac{K_1c_4\epsilon }{2}(\pi {m})^{1/2}}
\sum_{m'\in\mathbb{M}_n}(\pi {m'})^{4\gamma+1/2}
e^{-4b(\pi {m'})^{1/2}}\leq {C}.
\end{eqnarray*}
The term $II(m)$ is also bounded since $\epsilon$ is a constant.

\item[case $s>(1-s)_+\Leftrightarrow s>1/2$]
Here we choose $\epsilon$ such that 
$$4b(\pi {m''})^s-{K_1c_4\epsilon}(\pi {m''})^{(1-s)_+}/2=-4b(\pi {m''})^s$$ 
so that
\begin{eqnarray*}
  I(m)&\leq &c_3(\pi {m})^{4\gamma+\min(1-s,2-2s)}e^{-4b(\pi {m})^s}
\sum_{m'\in\mathbb{M}_n}e^{-\frac{K_1c_4\epsilon }{2}(\pi {m'})^{(1-s)_+}}\\ &&
+c_3e^{-\frac{K_1c_4\epsilon }{2}(\pi {m})^{(1-s)_+}}
\sum_{m'\in\mathbb{M}_n}(\pi {m'})^{4\gamma+\min(1-s,2-2s)}e^{-4b(\pi {m'})^s}\leq C.
\end{eqnarray*}
Moreover 
$$II(m)=\sum_{m'\in\mathbb{M}_n}\frac{1}{nC^2(\epsilon)}
e^{-K_2\sqrt{8b/K_1c_4}(\pi {m''})^{[s-(1-s)_+]/2}\sqrt{n}}\leq C.$$
\end{description}

In any case
$\epsilon =[8b/K_1c_4](\pi {m''})^{[s-(1-s)_+]_+}$, so that  
$$P_1(m,m')=C(q)(\pi {m''})^{[s-(1-s)_+]_+}\Delta^2(m'')/n$$
where $C(q)$ is a constant depending only on $k_0,k_1, b, \gamma, s$. 

\subsection{Proof of Proposition \ref{talagrand2}}
We split $Z_{n,2}(T)$ into two terms :
$$Z_{n,2}(T)=\frac{1}{2}Z_{n,2}^o(T)+\frac{1}{2}Z_{n,2}^e(T)$$
with $$\begin{cases}
\displaystyle Z_{n,2}^o(T)=\frac{2}{n}\sum_{i=1, i \text{ odd}}^{n}\left\{
T(X_i,X_{i+1})-\int T(x,y) F(x,y)dxdy\right\},\\
\displaystyle Z_{n,2}^e(T)=\frac{2}{n}\sum_{i=1, i \text{ even }}^{n}\left\{
T(X_i,X_{i+1})-\int T(x,y) F(x,y)dxdy\right\}.
\end{cases}$$

We bound $\E\left(\left[\sup_{T\in B(m')}(Z_{n,2}^{o})^{2}(T)-P_2(m,m')\right]_+
\1_{\Omega^*}\right)$. The second term can be bounded in the same way.
We write $Z_{n,2}^o(T)=(2/n)\sum_{i=1}^{n/2}\left\{T(V_i)-\E[T(V_i)]\right\}$
with $V_i=(X_{2i-1},X_{2i})$. 
In order to use  Lemma \ref{tala}, we introduce 
 $$ Z_{n,2}^{o*}(T)=\frac{1}{2}\nu_{n,1}(T)+\frac{1}{2}\nu_{n,2}(T) $$
where $$\begin{cases}
\displaystyle  \nu_{n,1}(T)=\frac{1}{p_n}\sum_{l=0}^{p_n-1}\frac{1}{q_n}
\sum_{i=2lq_n+1}^{(2l+1)q_n}\left\{T(V_i^*)-\E[T(V_i^*)]\right\},\\
\displaystyle  \nu_{n,2}(T)=\frac{1}{p_n}\sum_{l=0}^{p_n-1}\frac{1}{q_n}
\sum_{i=(2l+1)q_n+1}^{(2l+2)q_n}\left\{T(V_i^*)-\E[T(V_i^*)]\right\}.\\
\end{cases}$$
Since $X_i=X_i^*$ on $\Omega^*$, we can replace  $Z_{n,2}^{o}$ by $Z_{n,2}^{o*}$.
This leads us 
to bound $\E\left(\left[\sup_{T\in B(m')}\nu_{n,1}^{2}(T)-P_2(m,m')\right]_+\1_{\Omega^*}\right)$.
So we compute the bounds $M_1$, $H$ and $v$ of Lemma \ref{tala}.

\begin{enumerate}
\item If $T$ belongs to $\mathbb{S}_{m"}$,
$|T(x,y)|^2\leq \sum_{j,k}a_{j,k}^2\sum_{j,k}\varphi_{m'',j}^2(x)\varphi_{m'',k}^2(y)$
and so
$$\|T\|_\infty\leq \|T\|\|\sum_{j}\varphi_{m'',j}^2\|_\infty\leq \|T\|{m''},$$
using (1) of Lemma \ref{ordrevar}. Then 
$\|{1}/{q_n}\sum_{i=2lq_n+1}^{(2l+1)q_n}T\|_\infty\leq \|T\|{m''}$
and $M_1={m''}.$

\item Let us compute $H^2$.
  \begin{eqnarray*}
    \sup_{T\in B(m')}\nu_{n,1}^2(T)\leq \sum_{j,k}\nu_{n,1}^2(\varphi_{m'',j}\otimes \varphi_{m'',k})
  \end{eqnarray*}
Then, by taking the expectation,
  \begin{eqnarray*}
    \E\left(\sup_{T\in B(m')}\nu_{n,1}^2(T)\right)\leq \sum_{j,k}\frac{1}{p_n^2}
\var(\sum_{l=0}^{p_n-1}\frac{1}{q_n}
\sum_{i=2lq_n+1}^{(2l+1)q_n}\varphi_{m'',j}\otimes \varphi_{m'',k}(V_i^*))\\
\leq \sum_{j,k}\frac{1}{p_n^2}\sum_{l=0}^{p_n-1}\var(\frac{1}{q_n}
\sum_{i=2lq_n+1}^{(2l+1)q_n}\varphi_{m'',j}\otimes \varphi_{m'',k}(V_i^*)),
  \end{eqnarray*}
by using independence of the $A_l^*$.
Lemma \ref{mixing} then gives 
\begin{eqnarray*}
    \E\left(\sup_{T\in B(m')}\nu_{n,1}^2(T)\right)
\leq\frac{4}{p_nq_n}\|\sum_{j,k}(\varphi_{m'',j}\otimes \varphi_{m'',k})^2\|_\infty\sum{\beta_{2k}}
\leq\frac{16}{n}(\sum\beta_{2k}){m''}^2
  \end{eqnarray*}
Finally $H=4\sqrt{\sum\beta_{2k}}{{m''}}/{\sqrt{n}}$
\item $v$ remains to be calculated. If $T$ belongs to $B(m')$, using Lemma \ref{mixing}
\begin{eqnarray*}
&&\var[\frac{1}{q_n}\sum_{i=2lq_n+1}^{(2l+1)q_n}T(V_i^*)]\leq
\frac{4}{q_n}\E[T^2(V_1)b(V_1)]\\
&&\leq \frac{4}{q_n}\|T\|_\infty\sqrt{\E[T^2(V_1)]}\sqrt{\E[b^2(V_1)]}
\leq\frac{4}{q_n}\|T\|_\infty\sqrt{\|F\|_\infty}\sqrt{2\sum(k+1)\beta_{2k}}
\end{eqnarray*}
and so $v={4\|F\|_\infty^{1/2}}\sqrt{2\sum(k+1)\beta_{2k}}{m''}/q_n.$
\end{enumerate}

By writing Talagrand's inequality (Lemma \ref{tala}) with $\epsilon=1$, we obtain
$$\E\left([\sup_{T\in B(m')}(\nu_{n,1})^2(T)-6\frac{16}{n}(\sum\beta_{2k}){m''}^2]_+\1_{\Omega^*}\right)
\leq \frac{K}{n}\left\{{m''}e^{-K_1{m''}}+\frac{{m''}q_n^2}{n}e^{-K_2\sqrt{n}/q_n}\right\}.$$
Then by summation over $m'$
\begin{eqnarray*}
  \sum_{m'\in\mathbb{M}_n}\E\left([\sup_{T\in B(m')}(\nu_{n,1})^2(T)-\frac{96}{n}(\sum\beta_{2k}){m''}^2]_+
\1_{\Omega^*}\right)\\\leq \frac{K}{n}\{\sum_{m'\in\mathbb{M}_n}{m''}e^{-K_1{m''}}
+\sum_{m'\in\mathbb{M}_n}{m''}n^{2c-1}e^{-K_2n^{1/2-c}}\}\leq \frac{C}{n}
\end{eqnarray*}
since $c<1/2$. In the same way, we obtain 
\begin{eqnarray*}
  \sum_{m'\in\mathbb{M}_n}\E\left([\sup_{T\in B(m')}(\nu_{n,2})^2(T)-\frac{96}{n}(\sum\beta_{2k}){m''}^2]_+
\1_{\Omega^*}\right)\leq \frac{C}{n},
\end{eqnarray*}
which yields 
 $$\sum_{m'\in\mathbb{M}_n}\E\left([\sup_{T\in B(m')}(Z_{n,2}^o)^2(T)-P_2(m,m')]_+
\1_{\Omega^*}\right)\leq \frac{C}{n}$$
with $P_2(m,m')=96(\sum\beta_{2k}){m''}^2/n.$

\subsection{Proof of Corollary  \ref{ratesF}}

Let us compute the bias term. Since $F_m^*=F^*\1_{[-\pi m,\pi m]^2}$,
\begin{eqnarray*}
 && \|F-F_m\|^2=\frac{1}{4\pi^2}\iint_{([-\pi m,\pi m]^2)^c}|F^*(u,v)|^2dudv\\
&&\leq \frac{1}{4\pi^2}\iint_{[-\pi m,\pi m]^c\times\mathbb{R}}|F^*(u,v)|^2dudv
+ \frac{1}{4\pi^2}\iint_{\mathbb{R}\times[-\pi m,\pi m]^c}|F^*(u,v)|^2dudv
\end{eqnarray*}
 But
 \begin{eqnarray*}
   \iint_{[-\pi m,\pi m]^c\times\mathbb{R}}|F^*(u,v)|^2dudv\leq 
L((\pi m)^2+1)^{-\Delta}e^{-2A(\pi m)^R}.
 \end{eqnarray*}

Thus $\|F-F_m\|^2=O((\pi m)^{-2\Delta}e^{-2A(\pi m)^R})$ and 
$$\E\|F-\tilde{F}\|^2\leq C'\inf_{m\in\mathbb{M}_n}\left\{(\pi m)^{-2\Delta}e^{-2A(\pi m)^R}
+(\pi m)^{[s-(1-s)_+]_++4\gamma+2-2s}\frac{e^{4b(\pi m)^s}}{n}\right\}+\frac{C}{n}.$$

Next the bias-variance trade-off is performed similarly to \cite{lacourdeconvo}.

\subsection{Proof of Theorem \ref{mainpi}}

Let $E_n=\{\|f-\tilde{f}\|_\infty\leq f_0/{2}\}$.
On $E_n$ and for $x\in B$, $\tilde{f}(x)=\tilde{f}(x)-f(x)+f(x)\geq f_0/2$.
Since $\tilde{F}$ belongs to $\mathbb{S}_{\hat{M}}$, using \eqref{connex}, 
$\|\tilde{F}\|_\infty\leq \hat{M}\|\tilde{F}\|$. But \eqref{dd} gives 
$\|\tilde{F}\|\leq \Delta(\hat{M})$ so that $\|\tilde{F}\|_\infty\leq \hat{M}\Delta(\hat{M})$.
Since $\hat{M}$ belongs to $\mathbb{M}_n$, $\Delta(\hat{M})\leq \sqrt{n}$ 
and Lemma \ref{del} gives $\hat{M}\leq \Delta(\hat{M})^{1/(2\gamma+1)}$ if $s=0$ or 
$\hat{M}\leq (\ln\Delta(\hat{M}))^{1/s}$ otherwise.
So, for $n$ large enough, $(2/f_0)\|\tilde{F}\|_\infty\leq n$ and 
$\tilde{\Pi}(x,y)={\tilde{F}(x,y)}/{\tilde{f}(x)} .$

For all $(x,y)\in B^2$,
\begin{eqnarray*}
|\tilde{\Pi}(x,y)-\Pi(x,y)|^2&\leq &\left|\frac{\tilde{F}(x,y)-\tilde{f}(x)\Pi(x,y)}
{\tilde{f}(x)}\right|^2\1_{E_n}+(|\tilde{\Pi}(x,y)|+|\Pi(x,y|)^2\1_{E_n^C}\\
&\leq &\frac{|\tilde{F}(x,y)-F(x,y)+\Pi(x,y)(f(x)-\tilde{f}(x))|^2}
{f_0^2/4}\\ &&+2(\|\tilde{\Pi}\|_\infty^2+|\Pi(x,y)|^2)1_{E_n^C}.
\end{eqnarray*}
Since $\int_B \Pi^2(x,y)dy \leq \|\Pi\|_{B,\infty}\int_B \Pi(x,y)dy\leq \|\Pi\|_{B,\infty}$ for all $x\in B$,
\begin{eqnarray*}
\E\|\Pi-\tilde{\Pi}\|_B^2&\leq &\frac{8}{f_0^2}[\E\|F-\tilde{F}\|^2
  +\|\Pi\|_{B,\infty}\E\|f-\tilde{f}\|^2]+2|B|(|B|n^2+\|\Pi\|_{B,\infty})P(E_n^C).
\end{eqnarray*}
We still have to prove that $P(E_n^C)\leq Cn^{-3}$. Given that 
$\|f-\tilde{f}\|_\infty \leq \|f-f_{\hat{m}}\|_\infty+\|f_{\hat{m}}-\hat{f}_{\hat{m}}\|_\infty$
we obtain $$P(E_n^c)\leq P(\|f-f_{\hat{m}}\|_\infty>f_0/4)+P(\|f_{\hat{m}}-\hat{f}_{\hat{m}}\|_\infty>f_0/4).$$
Let us prove now that if $f$ belongs to $\mathcal{A}_{\delta,r,a}(l)$ with $\delta>1/2$,
$\|f-f_{m}\|_\infty=O(m^{1/2-\delta}e^{-a(\pi L_m)^r}).$
Since $f_m^*=f^*\1_{[-\pi m, \pi m]}$ and using the inverse Fourier transform, 
$$\|f-f_m\|_\infty\leq\frac{1}{2\pi}\int_{|u|\geq \pi m}|f^*(u)|du.$$
Let $\alpha\in (1/4,\delta/2)$. By considering that function 
$x\mapsto (x^2+1)^{\delta/2-\alpha}e^{a|x|^r}$ is increasing and using the Schwarz inequality, we obtain
$$\|f-f_m\|_\infty\leq\frac{1}{2\pi}((\pi m)^{2}+1)^{-\delta/2+\alpha}e^{-a(\pi m)^r}\sqrt{l}
\sqrt{\int_{|u|\geq \pi m} (u^2+1)^{-2\alpha}du}.$$
But $\int_{|u|\geq \pi m} (u^2+1)^{-2\alpha}du\leq C(\pi m)^{1-4\alpha}$ and then 
$$\|f-f_m\|_\infty\leq\frac{\sqrt{Cl}}{2\pi}((\pi m)^{2}+1)^{-\delta/2+\alpha}e^{-a(\pi m)^r}
(\pi m)^{1/2-2\alpha}=O(m^{1/2-\delta}e^{-a(\pi m)^r}).$$
Thus, since ${\hat{m}}\geq \ln \ln n$, $\|f-f_{\hat{m}}\|_\infty\to 0$ 
and for $n$ large enough $P(\|f-f_{\hat{m}}\|_\infty>f_0/4)=0.$
Next $$P(\|f_{\hat{m}}-\hat{f}_{\hat{m}}\|_\infty>f_0/4)\leq P(\Omega^{*c})+
P\left(\|f_{\hat{m}}-\hat{f}_{\hat{m}}\|\1_{\Omega^*}>\frac{f_0}{4\sqrt{{\hat{m}}}}\right).$$
Since $c\theta>3$, $P(\Omega^{*c})\leq Mn^{1-c\theta}\leq Mn^{-2}$.
We still have to prove that 
$$P\left(\|f_{\hat{m}}-\hat{f}_{\hat{m}}\|\1_{\Omega^*}>
\frac{f_0}{4\sqrt{{\hat{m}}}}\right)\leq \frac{C}{n^2}.$$
First, we observe that 
$$\|f_{\hat{m}}-\hat{f}_{\hat{m}}\|^2=\sum_{j\in \mathbb{Z}}\left(\frac{1}{n}\sum_{i=1}^{n}
v_{\varphi_{\hat{m}j}}(Y_i)-\E[v_{\varphi_{\hat{m}j}}(Y_i)]\right)^2=\sup_{t\in B_{\hat{m}}}\nu_n^2(t)$$
where
$  \nu_n(t)=\frac{1}{n}\sum_{i=1}^{n}v_{t}(Y_i)-\E[v_{t}(Y_i)], B_m=\{t\in S_m, \|t\|\leq 1\}.$

Then $$P\left(\|f_{\hat{m}}-\hat{f}_{\hat{m}}\|\1_{\Omega^*}>\frac{f_0}{4\sqrt{{\hat{m}}}}\right)
=P\left(\sup_{t\in B_{\hat{m}}}|\nu_n(t)|\1_{\Omega^*}>\frac{f_0}{4\sqrt{{\hat{m}}}}\right).$$
As previously, we split $\nu_n(t)$ into two terms 
$$\nu_n(t)=\frac{1}{2p_n}\sum_{l=0}^{p_n-1}\frac{1}{q_n}\sum_{i=2lq_n+1}^{(2l+1)q_n}v_{t}(Y_i)-\E[v_{t}(Y_i)]
+\frac{1}{2p_n}\sum_{l=0}^{p_n-1}\frac{1}{q_n}\sum_{i=(2l+1)q_n+1}^{(2l+2)q_n}v_{t}(Y_i)-\E[v_{t}(Y_i)]$$
and it is sufficient to study
$$P\left(\sup_{t\in B_{\hat{m}}}\left|\frac{1}{p_n}\sum_{l=0}^{p_n}\frac{1}{q_n}\sum_{i=2lq_n+1}^{(2l+1)q_n}
v_{t}(Y_i^*)-\E[v_{t}(Y_i^*)]\right|>\frac{f_0}{4\sqrt{{\hat{m}}}}\right).$$
We use inequality \eqref{talapr} in proof of Lemma \ref{tala} with $\eta=1$
and $\lambda=\cfrac{f_0}{8\sqrt{{\hat{m}}}}$: 
\begin{eqnarray*}
  P\left(\sup_{t\in B_{\hat{m}}}\left|\frac{1}{p_n}\sum_{l=0}^{p_n}\frac{1}{q_n}\sum_{i=2lq_n+1}^{(2l+1)q_n}
v_{t}(Y_i^*)-\E[v_{t}(Y_i^*)]\right|>2H+\lambda\right)\\\leq \exp\left(-K{p_n}\min
\left(\frac{\lambda^2}{v},\frac{\lambda}{M_1}\right)\right)
\end{eqnarray*}
Here, we compute
$$M_1=\sqrt{\Delta(\hat{m})};\quad H^2=8\sum_k\beta_k\frac{\Delta(\hat{m})}{n}
;\quad v=4\sum_k\beta_k\frac{\Delta(\hat{m})}{q_n}.$$
Thus $$P\left(\sup_{t\in B_{\hat{m}}}|\nu_n(t)|>2H+\cfrac{f_0}{8\sqrt{{\hat{m}}}}\right)\leq 
2\exp\left(-K'\min\left(\frac{n}{\hat{m}\Delta(\hat{m})},\frac{p_n}{\sqrt{\hat{m}\Delta(\hat{m})}}\right)\right).$$
Now we use the assumption  $\forall m \quad m\Delta(m)\leq n/(\ln n)^2$.
For $n$ large enough,
$2H=4\sqrt{2\sum_k\beta_k}{\sqrt{\Delta(\hat{m})}}/{\sqrt{n}}\leq {f_0}/(8\sqrt{\hat{m}}).$
So $$P\left(\sup_{t\in B_{\hat{m}}}|\nu_n(t)|>\cfrac{f_0}{4\sqrt{{\hat{m}}}}\right)\leq 
2\exp\left(-K'\min\left((\ln n)^2,n^{1/2-c}\ln n\right)\right)\leq \frac{C}{n^3}.$$

\subsection{ Technical Lemmas}\label{techlem}

\begin{lem}\label{ordrevar}
 For each $m\in\mathcal{M}_n$
\begin{enumerate}
\item $\|\sum_{j}\varphi_{m,j}^2\|_\infty= m$
\item $\u_{\varphi_{m,j}}(x)={\sqrt{m}}/({2\pi})\int_{-\pi}^{\pi} e^{-ijv}e^{ixvm}[q^*(-vm)]^{-1}dv$
\item $\|\sum_{j}|\u_{\varphi_{m,j}}|^2\|_\infty= \Delta(m)$
\end{enumerate}
where $\Delta(m)$ is defined in \eqref{Deltam}.
\end{lem}

\emph{Proof of Lemma \ref{ordrevar}: }
First we remark that 
\begin{eqnarray*}
  \varphi_{m,j}^*(u)&=&\int e^{-ixu}  \sqrt{{m}}\varphi({m }x-j)dx\\
&=&\frac{1}{\sqrt{m}}e^{-iju/m}\int e^{-ixu/m}  \varphi(x)dx
=\frac{1}{\sqrt{m}}e^{-iju/m}\varphi^*(\frac{u}{m}).
\end{eqnarray*}
Thus using the inverse Fourier transform 
$$\varphi_{m,j}(x)=\frac{1}{2\pi}\int e^{iux}\frac{1}{\sqrt{m}}e^{-iju/m}\varphi^*(\frac{u}{m})du
=\frac{1}{2\pi}\int_{-\pi}^{\pi}e^{-ijv}\sqrt{m}e^{ixvm}dv.$$
The Parseval equality yields  $\sum_{j}\varphi_{m,j}^2(x)=1/2\pi\int_{-\pi}^{\pi}|\sqrt{m}e^{ixvm}|^2dv=m.$
The first point is proved. Now we compute $\u_{\varphi_{m,j}}(x)$
\begin{eqnarray*}
\u_{\varphi_{m,j}}(x)&=&\frac{1}{2\pi}\int e^{ixu} \frac{\varphi_{m,j}^*(u)}{q^*(-u)}du
=\frac{1}{2\pi}\int e^{ixu}\frac{1}{\sqrt{m}}e^{-iju/m}\varphi^*(\frac{u}{m}) \frac{du}{q^*(-u)}\\
&=&\frac{\sqrt{m}}{2\pi}\int e^{-ijv}e^{ixvm}\frac{\varphi^*(v)}{q^*(-vm)}dv.
\end{eqnarray*}
But $\varphi^*(v)=\1_{[-\pi,\pi]}(v)$ and thus the second point is proved.
Moreover $\u_{\varphi_{m,j}}(x)$ can be seen as a Fourier coefficient. Parseval's formula then gives  
 \begin{eqnarray*}
   \sum_{j}|\u_{\varphi_{m,j}}(x)|^2=\frac{1}{2\pi}\int_{-\pi}^{\pi}\left|\sqrt{m}e^{ixvm}\frac{1}{q^*(-vm)}\right|^2dv
=\frac{m}{2\pi}\int\left|q^*(-vm)\right|^{-2}dv.
 \end{eqnarray*}
Therefore $\| \sum_{j}|\u_{\varphi_{m,j}}|^2\|_\infty
=1/2\pi \int_{-\pi m}^{\pi m}|q^*(-u)|^{-2}du=\Delta(m).$

\begin{lem}\label{del}
If $q$ verifies $|q^*(x)|\geq k_0(x^2+1)^{-\gamma/2}\exp(-b|x|^s),$ then 
  \begin{enumerate}
\item $\Delta(m)\leq c_1(\pi m)^{2\gamma+1-s}e^{2b(\pi m)^s},$
\item $\Delta_2(m)\leq c_2(\pi m)^{4\gamma+1-s}e^{4b(\pi m)^s}.$
  \end{enumerate}
Moreover if $|q^*(x)|\leq k_1(x^2+1)^{-\gamma/2}\exp(-b|x|^s),$ then $\Delta(m)\geq c'_1(\pi m)^{2\gamma+1-s}e^{2b(\pi m)^s}.$ 
\end{lem}

The proof of this result is omitted. It is obtained by distinguishing the cases 
$s>2\gamma+1$ and $s\leq 2\gamma+1$ and with standard evaluations of integrals.

\begin{lem}\label{tala}
Let $T_1,\dots,T_n$ be independent random variables and $\nu_n(r)=
(1/n)\sum_{i=1}^n[r(T_i)-\E(r(T_i)]$, for $r$ belonging to a countable class 
$\mathcal{R}$ of measurable functions. Then, for $\epsilon>0$,
$$\E[\sup_{r\in\mathcal{R}}|\nu_n(r)|^2-2(1+2\epsilon)H^2]_+\leq C
\left(\frac{v}{n}e^{-K_1\epsilon\frac{nH^2}{v}}+\frac{M_1^2}{n^2C^2(\epsilon)}
e^{-K_2C(\epsilon)\sqrt{\epsilon}\frac{nH}{M_1}}\right) $$
with $K_1=1/6$, $K_2=1/(21\sqrt{2})$, $C(\epsilon)=\sqrt{1+\epsilon}-1$ and $C$ a universal constant and where
$$\sup_{r\in\mathcal{R}}\|r\|_\infty\leq M_1, \quad  \E\left(\sup_{r\in\mathcal{R}}|\nu_n(r)|
\right)\leq H,\quad \sup_{r\in\mathcal{R}}\frac{1}{n}\sum_{i=1}^n\var(r(T_i))\leq v.$$
\end{lem}
Usual density arguments allow to use this result with non-countable class of functions $\mathcal{R}$.

\emph{Proof of Lemma \ref{tala}: }
We apply the Talagrand concentration inequality given in \cite{kleinrio} to the functions 
$s^i(x)=r(x)-\E(r(T_i))$ and we obtain
$$P(\sup_{r\in\mathcal{R}}|\nu_n(r)|\geq H+\lambda)\leq \exp\left(-\frac{n\lambda^2}{2(v+4HM_1)+6M_1\lambda}\right).$$
Then we modify this inequality following \cite{birge&massart98} Corollary 2 p.354. It gives
\begin{equation}
  \label{talapr}
  P(\sup_{r\in\mathcal{R}}|\nu_n(r)|\geq (1+\eta)H+\lambda)\leq \exp\left(-\frac{n}{3}\min
\left(\frac{\lambda^2}{2v},\frac{\min(\eta,1)\lambda}{7M_1}\right)\right).
\end{equation}
To conclude we set $\eta=\sqrt{1+\epsilon}-1$ and we use the formula
$\E[X]_+=\int_{0}^\infty P(X\geq t)dt$ with $X=\sup_{r\in\mathcal{R}}|\nu_n(r)|^2-2(1+2\epsilon)H^2.$

\begin{lem} (\cite{viennet97})\label{mixing}
Let $(T_i)$ a strictly stationary process with $\beta$-mixing
coefficients $\beta_k$. Then there exists a function $b$ such that 
$$\E[b(T_1)]\leq \sum_k\beta_k\quad \text{ and }\quad   
\E[b^2(T_1)]\leq 2\sum_k(k+1)\beta_k$$
and for all function $\psi$ (such that $\E[\psi^2(T_1)]<\infty)$ and for all $N$
$$\var(\sum_{i=1}^N \psi(T_i))\leq 4N\E[\psi^2(T_1)b(T_1)].$$
\end{lem}

\bibliographystyle{elsart-harv}
\bibliography{biblio}

\end{document}